\documentclass[12pt]{article}
\usepackage{amsmath,amsthm,amssymb}
\usepackage{amsfonts}
\newcommand{\en}{\enspace}

\begin{document}

\large

$$
$$

\medskip

\begin{center}\bf
AN UPPER BOUND FOR THE NUMBER \\
OF EIGENVALUES OF NON-SELFADJOINT \\
SCHR\"ODINGER OPERATOR
\end{center}

\smallskip

\begin{center}
\bf S. A. Stepin
\end{center}

\bigskip
\medskip

{\bf Abstract.} Estimates for the total multiplicity of eigenvalues for Schr\"o\-dinger operator are established 
in the case of compactly supported or exponentially decreasing complex-valued potential.

\medskip


\bigskip

\bigskip

Schr\"odinger operator $\,-\Delta+V\,$ with complex potential represents an adequate and capacious perturbation-theoretic model within non-selfadjoint setting (see [1]). Relative compactness of perturbation $\,V\,$ guarantees conservation of continuous spectrum filling the semi-axis $\,\mathbb R_+.\,$ A simple condition of relative compactness (see [2]) is given by integrability of potential being raised to an appropriate power which provides that operator $\,VR_0(\lambda)\,$ belongs to the corresponding Schatten - von Neumann class where $\,R_0(\lambda)\,$ denotes the resolvent of free Laplacian $\,-\Delta.$

Investigation of discrete spectral component of operator $\,-\Delta+V\,$ proves to be (see [3]) of a considerable interest and specific difficulty in the case of complex $\,V.\,$ Contemporary results on its location (the review of up-to-date state of research see in [4]) reduce to evaluation of the norm of the so-called Birman-Schwinger operator $\,V^{1/2}R_0(\lambda)|V|^{1/2}.\,$ Note that the estimates of Lieb-Thirring type obtained in non-selfadjoint case (see [5]) for the momenta of eigenvalues hardly enable one to directly evaluate their total number. At the same time results obtained in [6] on distribution of eigenvalues located in the sector containing $\,\mathbb R_+\,$ provide a certain information on their possible accumulation rate to continuous spectrum.

In the present paper an upper bound for the total number of eigenvalues (with multiplicities taken into account) is achieved for Schr\"odinger operator with compactly supported (section 1) and exponentially decaying (section 2) potential. To this end Nevanlinna-Jensen formula is applied to an appropriate Fredholm determinant vanishing at the eigenvalues (cf. [7]). Note that the determinant mentioned above proves to be resolvent denominator for the iterated Lippmann-Schwinger equation and plays here the same role as Jost function in one-dimentional case (cf. [8]).

In what follows for the sake of simplicity we consider three-dimensional configuration space; all the results formulated below admit natural extension to the case of arbitrary dimension.

\bigskip

{\bf 1.} In our setting the resolvent $\,R_0(\lambda)=(-\Delta-\lambda I)^{-1}\,$ is represented by an integral operator with the kernel
$\,\,\displaystyle{\frac{\exp\big(i\sqrt{\lambda}\,|x-y|\big)}{4\pi |x-y|}}\,.\,\,$
To derive an upper bound for the total multiplicity of eigenvalues for the operator $\,-\Delta+V\,$ with compactly supported potential we make use of the following

\medskip

{\bf Lemma 1.} {\it Given $\,\lambda=k^2,\,k\in\mathbb C_+,\,$ the integral kernel
$$
G_{\lambda}(x,y)\en=\en\frac1{16\pi^2}\int\frac{e^{ik|x-z|}e^{ik|z-y|}}{|x-z||z-y|}\,V(z)\,d^{\,3}z
$$
of the operator $\,R_0(\lambda)VR_0(\lambda)\,$ satisfies the estimate
$$
\big|G_{\lambda}(x,y)\big|\en\leqslant\en \frac{2\,C}{\sqrt{1+4|k|}-1}
$$
with constant
$$
C=\max\bigg\{\frac1{8\pi^2}\max_{x}\int\frac{|V(y)|\,d^{\,3}y}{|x-y|},\,\frac{\|V\|_{\infty}}{8\pi}\,+\,\frac{\sqrt{2}}{16}\,\|\nabla V\|_{\infty}\big(
{\rm diam}\,(\,{\rm supp}V)+1\big)\bigg\}.
$$
}

\smallskip

For $\,a\geqslant 0\,$ introduce a notation
$$
f(a)\en:=\en\sum_{n=0}^{\infty}\frac{n^{n/2}}{n!}\,a^n\en\leqslant\en (1+a)\exp(2a^2)
$$
and note that $\,f^{-1}(2)>1/2\,$ while $\,\ln\big(2-f(1/2)\big)\geqslant -3.$

\bigskip

{\bf Theorem 1.} \, {\it Suppose that potential $\,V\in{\rm C}^1(\mathbb R^3)\,$ is compactly supported. Then discrete spectrum of operator $\,-\Delta+V\,$ is located in the disc $\,\,|\lambda|\leqslant R\,\,$ of radius
$$
R\en\leqslant\en (C\|V\|_1)^2\big(\,1+C\|V\|_1\big)^2\,.
$$
Given arbitrary $\,\,\varepsilon>0\,\,$ and
$\,\,T>\max\big\{2R/\varepsilon-\varepsilon/8,2C\|V\|_1(1+2C\|V\|_1)\big\}\,\,$ for the total multiplicity $\,\,\,N(V)\,\,$ of the eigenvalues of the operator $\,\,-\Delta+V\,\,$ the estimate
$$
N(V)\,\leqslant\,\bigg(\!\ln\frac{T+\varepsilon/4}{\sqrt{T^2+R}}\bigg)^{-1}\!\bigg[\,\ln f\bigg(\frac{AB}{2\pi\varepsilon}\bigg)\,-\,\,\ln\bigg\{2\,-f\bigg(\frac{2C\|V\|_1}{\sqrt{1+4T}-1}\bigg)\bigg\}\bigg],
$$
holds where
$$
A=A(\varepsilon)\,:=\,\max |V(x)|\,e^{\varepsilon |x|}\,,\quad B=B(\varepsilon)\,:=\,\int |V(x)|\,e^{\varepsilon |x|}\,d^{\,3}x\,.
$$
}

\medskip

The first assertion of the theorem is obtained by the usage of uniform boundedness property valid for the eigenfunctions of operator  $\,-\Delta+V\,$ which belong to $\,\ker\big(I-(R_0(\lambda)V)^2\big).\,$ Namely it reduces to norm estimation of operator
$\,K(\lambda)=(R_0(\lambda)V)^2\!:{\rm L}_{\infty}(\mathbb R^3)\to{\rm L}_{\infty}(\mathbb R^3)\,$ with regard to lemma 1.

Given $\,\lambda=k^2,\,k\in\mathbb C_+,\,$ due to the equality
$$
\int\int\bigg|\frac{e^{ik\,|x-y|}}{|\,x-y|}\,V(y)\bigg|^2\,d^{\,3}x\,d^{\,3}y\en=\en\frac{2\pi}{{\rm Im}\,k}\,\|V\|_2^2\,\,<\,\,\infty
$$
the operator $\,R_0(\lambda)V\,$ belongs to Hilbert-Schmidt class and therefore operator $\,K(\lambda)\,$ is of trace class. The second claim of theorem 1 is deduced from the estimate of the value $\,N(V)\,$ by the number of zeroes of Fredholm determinant
\begin{multline*}
D(\sqrt{\lambda})\en=\en\det\big(I-K(\lambda)\big)\en=\\
=\en 1\en+\en\sum_{n=1}^{\infty}\frac{(-1)^n}{n!}\int\ldots\int
\det\big(G_{\lambda}(x_i,x_j)V(x_j)\big)\,d^{\,3}x_1\ldots d^{\,3}x_n\,.
\end{multline*}
Analytic in $\,\mathbb C_+\,$ function $\,D(k)\,$ admits continuation to $\,\mathbb C_-,\,$ so that for $\,{\rm Im}\,k>-\varepsilon/4\,$
the inequality
$$
|D(k)|\en\leqslant\en f\bigg(\frac{AB}{2\pi\varepsilon}\bigg)\,
$$
is valid (see [9]).
Setting $\,\varphi(z)=D(z+iT)\,$ and choosing $\,\rho\in(\sqrt{T^2+R}, T+\varepsilon/4)\,$ by virtue of  Nevanlinna-Jensen formula one has
$$
N(V)\,\ln\frac{\rho}{\sqrt{T^2+R}}\en\leqslant\en\frac1{2\pi}\int_0^{2\pi}\!\ln|\varphi(\rho e^{i\theta})|\,d\theta\,\,-\,\,\ln|\,\varphi(0)|\,.
$$
Lower bound for $\,|\varphi(0)|=|D(iT)|\,$ is derived from lemma 1 by usage of Hadamard inequality (see [10]) according to which
$$
|D(k)-1|\en\leqslant\en f\bigg(\frac{2C\|V\|_1}{\sqrt{1+4|k|}-1}\bigg)\,-\,1\,,\quad k\in\mathbb C_+\,.
$$

\bigskip

{\bf Corollary 1.} {\it For arbitrary $\,\varepsilon>0\,$ the inequality
\begin{multline*}
N(V)\,\leqslant\,2\bigg(\!\ln\bigg(1\,+\,\frac{(1+C\|V\|_1)^{-2}}4\min\bigg\{\frac{\varepsilon}{2C\|V\|_1},
\frac{(1+C\|V\|_1)^2}{1+2C\|V\|_1}\bigg\}^2\bigg)\bigg)^{-1}\\
\times\,\bigg[\,\,3\,+\,2\bigg(\frac{AB}{2\pi\varepsilon}\bigg)^2+\,\,\ln\bigg\{1\,+\,\,\frac{AB}{2\pi\varepsilon}\bigg\}\,\bigg]
\end{multline*}
holds true provided that hypothesis of theorem {\rm 1} is satisfied.
}

\medskip

Indeed assuming $\,\varepsilon/2\leqslant L:=C\|V\|_1(1+C\|V\|_1)\,$ set $\,T=4L^2/\varepsilon\,$ so that
$$
\frac{T+\varepsilon/4}{\sqrt{T^2+L^2}}\en=\en\bigg(1+\frac{\varepsilon^2}{16L^2}\bigg)^{1/2}\!;
$$
provided that $\,\varepsilon/2\geqslant L\,$ choose $\,T=2C\|V\|_1(1+2C\|V\|_1)\,$ and therefore
$$
\frac{T+\varepsilon/4}{\sqrt{T^2+L^2}}\en\geqslant\en\bigg(1+\frac{(1+C\|V\|_1)^2}{4(1+2C\|V\|_1)^2}\bigg)^{1/2}\!.
$$
However in both cases one has
$$
\ln\bigg(2-f\bigg(\frac{2C\|V\|_1}{\sqrt{1+4T}-1}\bigg)\bigg)\,\geqslant\,
\,\ln\big(2-f(1/2)\big)\,\geqslant\,-3\,.
$$

\medskip

{\bf 2.} The analogues of lemma 1 and theorem 1 (as well as its corollary) are formulated below for the case when $\,|V(x)|\leqslant A e^{-\varepsilon|x|}\,$ with certain $\,A>0\,$ and $\,\varepsilon>0.\,$
Denote by $\,h_{\varepsilon}\,$ the function inverse to $\,g_{\varepsilon}(t)=t\exp \varepsilon t\,$ so that  $\,h_{\varepsilon}(t)\leqslant t\,$ and besides $\,h_{\varepsilon}(t)\sim\,\varepsilon^{-1}\ln t\,$ when $\,t\to\infty.$
\bigskip

{\bf Lemma 2.} \, {\it For $\,\lambda=k^2,\,k\in\mathbb C_+,\,$ integral kernel
$\,G_{\lambda}(x,y)\,$
of operator $\,R_0(\lambda)VR_0(\lambda)\,$ admits the estimate
$$
|G_{\lambda}(x,y)|\en\leqslant\en \frac{\widetilde{C}}{h_{\varepsilon}(|k|)}
$$
with constant
$$
\widetilde{C}\,=\, \frac1{8\pi^2}\,\max\bigg\{\max_{x}\int|V(y)|\,\frac{d^{\,3}y}{|x-y|}\,,\,\,\pi A\,\big(1+\sqrt{2}\pi\big)\bigg\}\,.
$$
}

\bigskip

Proof of lemma 2 as well as that of lemma 1 make use of the inequality
$$
|\,G_{\lambda}(x,y)|\en\leqslant\en\frac1{4\pi^2|x-y|}\,\max_{x}\int|V(y)|\,\frac{d^{\,3}y}{|x-y|}
$$
and the following (cf. [9]) general statement giving the estimate for the integral kernel of operator $\,R_0(\lambda)VR_0(\lambda)\,$ applied to compactly supported or exponentially decreasing potential $\,V\,$  respectively.

\bigskip

{\bf Proposition.} {\it Suppose that potential $\,V\in{\rm C}^1(\mathbb R^3)\,$ is bounded. Then given arbitrary $\,k\in\mathbb C_+\,$ for the integral kernel $\,G_{\lambda}(x,y),$ $\lambda=k^2,\,$ the inequality
$$
|\,G_{\lambda}(x,y)|\en\leqslant\en\frac1{8|k|}\bigg(\frac1{\pi}\,\|V\|_{\infty}\,+\,\,\frac1{\sqrt{2}}\int_c^{\infty}\!\!\max_{E(x,y,r)}\!|\nabla V|\,\frac{r\,dr}{\sqrt{r^2-c^2}}\bigg)\,,
$$
is valid where $\,c=|x-y|/2\,$ while $\,E(x,y,r)\,$ is an ellipsoid of revolution with foci $\,x\,$ and $\,y\,$ and semiaxes $\,r\,$ and $\,\sqrt{r^2-c^2}.$
}

\bigskip

To derive lemma 1 from the proposition formulated above set $\,d=\,{\rm diam}({\rm supp}V),$ $\,\delta=\,{\rm dist}([x,y],{\rm supp}V)\,$ and take into account that $\,V\,$ vanishes on ellipsoids $\,E(x,y,r)\,$ when $\,r<\sqrt{c^2+\delta^2}\,$ or $\,r>c+d+\delta\,$ and hence
$$
\int_c^{\infty}\!\!\max_{E(x,y,r)}\!|\nabla V|\,\frac{r\,dr}{\sqrt{r^2-c^2}}\en\leqslant\en\|\nabla V\|_{\infty}(c+d)\,.
$$

In order to prove lemma 2 let $\,\mu=\min\limits_{t\in[0,1]}|\,x+t(y-x)|\,$ and note that the following, involving standard Heaviside step function $\,\vartheta,\,$ estimate
$$
|\nabla V|\en\leqslant\en\varepsilon A\,\Big(\,e^{\varepsilon(c+\mu-r)}\vartheta(r-\sqrt{\mu^2+c^2})\,       +\,e^{\varepsilon(\sqrt{r^2-c^2}-\mu)}\vartheta(c+\mu-r)\Big)
$$
valid on ellipsoid $\,E(x,y,r)\,$ implies the inequality
$$
\int_c^{\infty}\!\!\max_{E(x,y,r)}\!|\nabla V|\,\frac{r\,dr}{\sqrt{r^2-c^2}}\en\leqslant\en A\,e^{\varepsilon c}\,.
$$

\bigskip

{\bf Theorem 2.} \, {\it Let $\,V\,$ be continuously differentiable potential satisfying condition
$$
\big|\nabla V(x)\,\big|\en\leqslant\en \varepsilon A\,e^{-\varepsilon |x|}
$$
with a certain constant $\,A>0\,$ and moreover $\,V(x)\to 0\,$ at infinity.
Then discrete spectrum of operator $\,-\Delta+V\,$ is located in a disc $\,\,|\lambda|\leqslant R\,\,$ of radius
$\,\,R\,\leqslant\,\big(\widetilde{C}\|V\|_1\big)^2\exp\big(2\,\varepsilon\,\widetilde{C}\|V\|_1\big).\,\,$
If additionally
$$
B\,=\,\int |V(x)|\,e^{\varepsilon |x|}\,d^{\,3}x\,<\,\infty\,
$$
then for arbitrary $\,\,T>\max\big\{2R/\varepsilon-\varepsilon/8,\,g_{\varepsilon}\big(2\widetilde{C}\|V\|_1\big)\big\}\,\,$ and  $\,\,\rho\in[\sqrt{T^2+R},$ $  T+\varepsilon/4]\,$ the total multiplicity $\,N(V)\,$ of eigenvalues of operator $\,-\Delta+V\,$ admits the following estimate
$$
N(V)\en\leqslant\en\bigg(\ln\frac{\rho}{\sqrt{T^2+R}}\bigg)^{-1}\bigg[\,\ln f\bigg(\frac{AB}{2\pi\varepsilon}\bigg)\,
-\,\,\ln\bigg\{2\,-f\bigg(\frac{\widetilde{C}\|V\|_1}{h_{\varepsilon}(T)}\bigg)\bigg\}\bigg].
$$
}

\bigskip

Under the assumptions of theorem 2 the inequality
$$
|V(x)|\en\leqslant\en A\,e^{-\varepsilon|x|}\,,
$$
holds which enables one (cf. sketch of the proof of theorem 1) to carry out analytic continuation of the corresponding Fredholm determinant $\,D(k)\,$ from $\,\mathbb C_+\,$ into the half-plane $\,{\rm Im}\,k>-\varepsilon/4\,$ with the validity of boundedness condition
$$
|D(k)|\en\leqslant\en f\bigg(\frac{AB}{2\pi\varepsilon}\bigg)\,.
$$
Note that for super-exponentially decreasing potentials $\,V(x)\,$ this approach can be further applied to sharpen the upper bound for $\,N(V)\,$ obtained in theorem 2.

\bigskip
\medskip

{\bf Corollary 2.} {\it The estimate 
\begin{multline*}
N(V)\en\leqslant\en2\,\bigg(\ln\bigg(1\,+\,\frac14\min\bigg\{1\,,\,\frac{\varepsilon}{2\widetilde{C}\|V\|_1}\bigg\}^2\exp\big(\!-\!2\widetilde{C}\|V\|_1\big)\bigg)\bigg)^{-1}\times\\
\times\,\bigg[\,\,3\,+\,2\bigg(\frac{AB}{2\pi\varepsilon}\bigg)^2+\,\,\ln\bigg\{1\,+\,\,\frac{AB}{2\pi\varepsilon}\bigg\}\,\bigg]
\end{multline*}
is valid provided that hypothesis of theorem {\rm 2} is satisfied.  
}

\bigskip
\medskip

Indeed assuming $\,\varepsilon\leqslant 2\widetilde{C}\|V\|_1\,$ set $\,\rho=T+\varepsilon/4\,$ where $\,T=4M^2/\varepsilon\,$ and  $\,M=g_{\varepsilon}(\widetilde{C}\|V\|_1)\,$ so that
$$
\frac{\rho}{\sqrt{T^2+M^2}}\en=\en\bigg(1+\frac{\varepsilon^2}{16M^2}\bigg)^{1/2}\,\,,\quad f\bigg(\frac{\widetilde{C}\|V\|_1}{h_{\varepsilon}(T)}\bigg)\,\leqslant\, f\big(1/2\big)\,.
$$

\medskip

\noindent In case $\,\varepsilon\geqslant 2\widetilde{C}\|V\|_1\,$ choose $\,T=g_{\varepsilon}(2\widetilde{C}\|V\|_1)\,$  
let $\,\rho=T+\varepsilon/4\,$ and verify
$$
\frac{\rho}{\sqrt{T^2+M^2}}\en\geqslant\en\bigg(1+\frac14\exp\big(-\widetilde{C}\|V\|_1\big)\bigg)^{1/2}\,,\quad f\bigg(\frac{C\|V\|_1}{h_{\varepsilon}(T)}\bigg)\,=\, f\big(1/2\big)\,.
$$

\bigskip
\bigskip

{\bf Remark. } Within the setting in question for the eigenvalues $\,\lambda\in{\mathbb C}\setminus\mathbb R_+\,$ of operator $\,-\Delta+V\,$ an upper bound
$$
|\lambda|\en\leqslant\en D(\gamma)\,\bigg(\int |V(x)|^{\gamma+3/2}d^{\,3}x\bigg)^{1/\gamma}
$$
is established in [4] involving a certain constant  $\,D(\gamma),\,\gamma\in(0,1/2].\,$ The estimates for the "discrete spectral radius" $\,R\,$ obtained in theorems 1 and 2 supplement and somewhat (e.g. for  $\,V\,$ small enough) strengthen the bound mentioned here. As regards homogeneity degree with respect to potential the evaluation of $\,R\,$ from theorem 1 for large $\,V\,$ is consistent with the best one among these indicated above (corresponding to the exponent $\,\gamma=1/2$).

\bigskip
\bigskip
\newpage

\begin{center}
\bf REFERENCES
\end{center}

\smallskip

\begin{enumerate}
\item
I.M.Gelfand \, On the spectrum of non-selfadjoint operators, Russian Math. Surveys, 1952. V.7. N6. P.183-184.
\item
E.H. Lieb, W. Thirring \, Inequalities for the moments of the eigenvalues of the Schr\"odinger Hamiltonian and their relation to Sobolev inequalities, Studies in Mathematical Physics, NJ: Princeton Univ. Press, 1976. P.269-303.
\item
M. Demuth, M. Hansmann, G. Katriel \, On the discrete spectrum of non-selfadjoint operators, J. Funct. Anal., 2009. V.257. N9. P.2742-2759.
\item
R.L. Frank \, Eigenvalue bounds for Schr\"odinger operators with complex potentials, Bull. London Math. Soc., 2011. V.43. N4. P.745-750.
\item
R.L. Frank, A. Laptev, E.H. Lieb, R. Seiringer \, Lieb-Thirring inequalities for Schr\"odinger operators with complex-valued potentials, Lett. Math. Phys., 2006. V.77. N3. P.309-316.
\item
A. Laptev, O. Safronov \, Eigenvalue estimates for Schr\"odinger operators with complex potentials, Comm. Math. Phys., 2009. V.292. N1. P.29-54.
\item
S.A. Stepin \, On spectral components of the Schr\"odinger operator with complex potential, Russian Math. Surveys, 2013. V.68. N1. P.186-188.
\item
B. Simon \, Resonances in one dimension and Fredholm determinants, J. Funct. Anal., 2000. V.178. N2. P.396-420.
\item
R.M. Martirosyan \, On the spectrum of various perturbations of the Laplace operator in spaces of three and more dimensions, Izv. Akad. Nauk, 1960. V.24. N6. P.897-920.
\item
W.V. Lovitt \, Linear integral equations, NY: Mc Graw-Hill, 1924.
\end{enumerate}

\end{document}